\documentclass[10pt]{amsart}
\usepackage{amsmath,amssymb,enumerate}

\usepackage{epsfig,fancyhdr,color}

\newtheorem{theorem}{\rm\bf Theorem}[section]
\newtheorem{proposition}[theorem]{\rm\bf Proposition}
\newtheorem{lemma}[theorem]{\rm\bf Lemma}
\newtheorem{corollary}[theorem]{\rm\bf Corollary}

\theoremstyle{definition}
\newtheorem{definition}[theorem]{\rm\bf Definition}

\theoremstyle{remark}

\newtheorem{example}[theorem]{\rm\bf Example}

\def\interieur#1{\mathord{\mathop{\kern 0pt #1}\limits^\circ}}

\title[A rigidity theorem for the mapping class group]{A rigidity theorem for the mapping class group action on the space of unmeasured foliations on a surface}
 
\author{Athanase Papadopoulos}
\address{A. Papadopoulos, Institut de Recherche Math{\'e}matique Avanc\'ee,
Universit{\'e} Louis Pasteur and CNRS,
7 rue Ren\'e Descartes,
 67084 Strasbourg Cedex - France} \email{papadopoulos@math.u-strasbg.fr}

\date{\today}

\begin{document}

\begin{abstract} Let $S$ be a surface of finite type   which is not a sphere with at most four punctures, a torus with at most two punctures, or a closed surface of genus two. Let $\mathcal{MF}$ be the space of equivalence classes of measured foliations of compact support on $S$ and let  $\mathcal{UMF}$ be the quotient space of $\mathcal{MF}$  obtained by identifying two equivalence classes whenever they can be represented by topologically equivalent foliations, that is, forgetting the transverse measure. The extended mapping class group $\Gamma^*$ of $S$ acts as by homeomorphisms of $\mathcal{UMF}$. We show that the restriction of the action of the whole homeomorphism group of  $\mathcal{UMF}$ on some dense subset of  $\mathcal{UMF}$ coincides with the action of $\Gamma^*$ on  that subset. More precisely, let  $\mathcal{D}$ be the natural image in $\mathcal{UMF}$ of the set of homotopy classes of not necessarily connected essential disjoint and pairwise nonhomotopic simple closed curves on $S$. The set $\mathcal{D}$ is dense in $\mathcal{UMF}$, it is invariant by the action of $\Gamma^*$ on $\mathcal{UMF}$ and  the restriction of the action of $\Gamma^*$  on $\mathcal{D}$ is faithful. We prove that 
 the restriction of the action on $\mathcal{D}$ of the group $\mathrm{Homeo}(\mathcal{UMF})$ coincides with the action of $\Gamma^*(S)$ on that subspace.
\bigskip


\end{abstract}

\maketitle

\section{introduction}\label{s1}

In this paper, $S=S_{g,p}$ is an oriented  surface of finite type, of genus $g\geq 0$ with $p\geq 0$ punctures. 
A simple closed curve on $S$ is said to be {\it essential} if it  is not homotopic to a point or to a puncture of $S$. We denote by
$\mathcal{S}$ the set of isotopy classes of essential simple closed curves on $S$ and by
 $\mathcal{MF}=\mathcal{MF}(S)$ the space of compactly supported measured foliations on $S$ up to the equivalence relation  $\sim$ generated by isotopy and Whitehead moves.  $\mathbb{R}_+$ denotes the set of nonnegative reals and $\mathbb{R}^{\mathcal{S}}_+$ is  the set of all functions from $\mathcal{S}$ to $\mathbb{R}_+$. We recall that the space $\mathcal{MF}$ is
equipped with a topology, defined by Thurston,  which is induced  from the embedding $\mathcal{MF}\to \mathbb{R}_+^{\mathcal{S}}$ by means of the intersection function $i:\mathcal{MF}\times\mathcal{S}\to\mathbb{R}_+$  (see \cite{Thurston1988} and \cite{FLP}).
  We also recall that the set $\mathcal{S}$ admits a natural embedding in the space $\mathcal{MF}$, and that the intersection function $i:\mathcal{MF}\times\mathcal{S}\to\mathbb{R}_+$  extends to an intersection function defined on $\mathcal{MF}\times\mathcal{MF}$, which we denote by the same letter $i:\mathcal{MF}\times\mathcal{MF}\to\mathbb{R}_+$. The quotient of $\mathcal{MF}$ by the natural action of the group $\mathbb{R}_+^*$ of positive reals is the space $\mathcal{PMF}=\mathcal{PMF}(S)$ of projective classes of compactly supported measured foliations on $S$.

The {\it mapping class group} $\Gamma(S)$ of $S$ is the group of isotopy classes of orientation-preserving homeomorphisms of $S$. The {\it extended mapping class group} $\Gamma^*(S)$ of $S$ is the group of all isotopy classes of homeomorphisms of $S$.

Let $\mathcal{UMF}=\mathcal{UMF}(S)$ be the quotient of $\mathcal{MF}(S)$ obtained by identifying two elements of $\mathcal{UMF}$ whenever these elements can be represented by topologically equivalent foliations (that is, the foliations are the same if we forget about the transverse measures).  The extended mapping class group $\Gamma^*(S)$ admits a natural action by homeomorphisms on
$\mathcal{UMF}$. This action is induced by the one which consists of taking images of foliations on $S$ by homeomorphisms of $S$.
 I started several years ago a study of the dynamics of the action of the mapping class group on $\mathcal{UMF}$ (see \cite{Papadopoulos1995}).  It is easy to see that this space, equipped with the quotient of the topology of $\mathcal{MF}$, is non-Hausdorff. We shall exploit this non-Hausdorffness in the proof of the theorem below.

We denote by $\mathrm{Homeo}(\mathcal{UMF})$ the group of homeomorphisms of $\mathcal{UMF}$.

In this paper, we prove the following result, which says that in some sense the action on $\mathcal{UMF}$ of the group $\mathrm{Homeo}(\mathcal{UMF})$ coincides on a dense subset of $\mathcal{UMF}$ with the action of the extended mapping class group $\Gamma^*(S)$ on that space.

\begin{theorem}\label{th} Let $S$ be a surface which is not a sphere with at most four punctures or a torus with at most two punctures. Then, there exists a dense subset $\mathcal{D}$ of $\mathcal{UMF}$ which is invariant by the group $\mathrm{Homeo}(\mathcal{UMF})$ and which has the following properties:
\begin{enumerate}
\item \label{i1} If $h$ is any homeomorphism of $\mathcal{UMF}$, then there exists an element $h^*$ of $\Gamma^*(S)$ such that the restriction on $\mathcal{D}$ of the actions of $h$ and $h^*$ coincide; 
\item  \label{i2} Suppose furthermore that $S$ is not a closed surface of genus 2. Then, if$h_1$ and $h_2$ are distinct elements of $\Gamma^*(S)$, their induced actions on $\mathcal{D}$ are different. In particular, the natural homomorphism from $\Gamma^*(S)$  to $\mathrm{Homeo}(\mathcal{UMF})$ is injective.
\end{enumerate}
\end{theorem}

 We shall see that the set $\mathcal{D}$ in the statement consists of the natural image in $\mathcal{UMF}$ of the set of systems of curves on $S$, that is, finite collections of pairwise non-isotopic and disjoint essential simple closed curves.  

The proof of Theorem \ref{th} is given in the last section of this paper. It uses the action of $\Gamma^*(S)$ on the curve complex of $S$. We recall that
the curve complex $C(S)$ of $S$ is the abstract simplicial complex on the vertex set $\mathcal{S}$ whose $k$-simplices, for all $k\geq 0$, are the subsets of $\mathcal{MF}$ that can be represented  by $k+1$ disjoint pairwise non-homotopic essential simple closed curves on $S$. The complex $C(S)$ has been introduced by Harvey in \cite{Harvey}. Masur and Minsky showed that $C(S)$,  equipped with its natural simplicial metric, is Gromov hyperbolic, and 
 Klarreich identified the Gromov boundary of
$C(S)$ as the subspace of $\mathcal{UMF}$ consisting of equivalence classes of minimal foliations (that is, foliations in which every leaf, including the singular leaves, is dense).

To prove property (\ref{i1}) in the statement Theorem \ref{th}, we first prove that any homeomorphism $h$ of $\mathcal{UMF}$ preserves the natural image in that space of set $\mathcal{S}'$ of isotopy classes of systems of curves on $S$, a system of curves being a collection of disjoint pairwise non-isotopic essential simple closed curves on $S$. This is done by 
defining a notion of {\it adherence number} for the points in $\mathcal{UMF}$, which is invariant by homeomorphism of $\mathcal{UMF}$ and which is a measure for non-Hausdorffness at these points. The elements of $\mathcal{UMF}$ which represent 
elements in $\mathcal{S}'$ are those that have the maximal adherence number. Then, we show that for each $k\geq 1$, $h$ preserves the natural image in $\mathcal{UMF}$  of the set $\mathcal{S}_k$ of isotopy classes of systems of curves on $S$ which have $k$components. This is done by considering the notion of {\it adherence set} for the  elements in $\mathcal{UMF}$, and showing that for $F\in\mathcal{S}_k$ and $F'\in\mathcal{S}'$ with $k\not= k'$, the adherence sets of $F$ and $F'$ are not homeomorphic.  Once we know that $h$ acts on each subset $\mathcal{S}_k$ of $\mathcal{UMF}$, we can define an action of $h$ on the curve complex $C(S)$ of $S$, and we apply a theorem of Ivanov, Korkmaz and Luo saying that except for the surfaces excluded in the hypothesis of Theorem  \ref{th}, each automorphism of the curve complex is induced by an element of the extended mapping class group. Property (\ref{i2})  in the  theorem is easy to prove.

\section{Adherence in topological spaces}

In this section, $X$ is a topological space. We introduce the notion  of adherent points, of adherence of a set and of adherence number of a point in $X$. These notions are interesting only in the case where $X$ is not Hausdorff. In some sense, the adherence number at a point measures a degree of non-Hausdorffness of $X$ at that point.  These notions will be essential in the proof of Theorem  \ref{th} that we give below.

\begin{definition}[Adherence]
Let $x$ and $y$ be two points in $X$.
We say that $x$ is {\it adherent to} $y$ in $X$ if every neighborhood of $x$ intersects every neighborhood of $y$.
\end{definition}

\begin{definition}[Adherence set]
Let $x$ be a point in $X$.
The {\it adherence set} $\mathcal{A}(x)$ of $x$ is the set of elements in $X$ which are adherent to $x$.
\end{definition}

\begin{definition}[Complete adherence set]
A subset $Y$ of $X$ is a {\it complete adherence set} in $X$ if for any two elements $x$ and $y$ of $Y$, $x$ is adherent to $y$ in $X$. \end{definition}

\begin{definition} Let $x$ be a point in $X$. The
{\it adherence number} $\mathcal{N}(x)$ of $x$ is the element of $\mathbb{N}\cup\{\infty\}$ defined by
\[\mathcal{N}(x)=\sup\{ \mathrm{Card}(A)\ \vert \ x\in A\text{ and } A \text{ is an adherence set in $X$}\}.\]
\end{definition}

In the next section, we shall deal with all these notions in the setting of the non-Hausdorff space $\mathcal{UMF}$. Let us start by mentioning a simple concrete example.

\begin{example}
A standard example of a non-Hausdorff space is the ``real line with two origins", obtained by considering two copies $R_1\simeq \mathbb{R}$ and $R_2\simeq \mathbb{R}$ of the real line $\mathbb{R}$, equipped with the natural bijection $f:R_1\to R_2$, and taking $X$ to be the quotient of the disjoint union $R_1\cup R_2$ by the map which identifies each point in $R_1$ that is distinct from the origin  with its image $f(x)$ in $R_2$. The space $X$ is non-Hausdorff. Let $\pi: R_1\cup R_2\to X$ be the natural map. If $x\in X$ is not the image by $\pi$ of the origin of $R_1$ or of $R_2$, then $\mathcal{N}(x)=1$. In the remaining case, $\mathcal{N}(x)=2$.
\end{example}

\section{Adherence in $\mathcal{UMF}$}

We need to recall the following technical detail that concerns measured foliations on punctured surfaces. We say that a measured foliation $F$ on $S$ has {\it compact support} if each puncture of $S$ has an annular neighborhood on which the foliation induced by $F$ is a foliation by closed leaves which are all parallel to that puncture. 
We note that this condition means that if we equip the surface with a complete finite volume hyperbolic structure and if we replace each leaf of $F$ by its geodesic representative when such a representative exists and with the empty set if such a representative does not exist, then the geodesic lamination that we obtain is compactly supported in the usual sense. Note also that we do not include in our space $\mathcal{UMF}$ the equivalence class of the measured foliation whose leaves are all parallel to punctures of $S$, that is, we do not include the ``empty" measured foliation (or lamination) with compact support. 

We also note that the isotopies and Whitehead moves that generate the equivalence relation $\sim$ must preserve transverse measure, except for the foliated annuli that are neighborhoods of punctures, on which the transverse measure is irrelevant.

It is convenient to represent the elements of $\mathcal{MF}$ by {\it partial} measured foliations, that is, measured foliations whose supports are nonempty (and not necessarily connected) subsurfaces with boundary of $S$. 
We shall denote by $\mathrm{Supp}(F)$ the support of a partial measured foliation $F$.
If $F$ is a partial measured foliation, then, a genuine measured foliations $F_0$ (which we sometimes call a {\it total} measured foliation, to stress the fact that its support is equal to $S
$) is obtained from $F$ by first inserting around each puncture an annulus foliated by closed curves parallel to that puncture, assigning to that annulus an arbitrary transverse measure, and then collapsing each connected component of $S\setminus F$ onto a spine. The equivalence class of $F_0$ does not depend on the choice of these spines, and in this way a partial measured foliation gives a well-defined element of $\mathcal{MF}$.

Given a (partial) measured foliations $F$ on $S$, we shall use the notation $[F]$ for its equivalence class in both spaces $\mathcal{MF}$ and  $\mathcal{UMF}$, and we shall make the ambient space  clear when using this notation.

 Let $F$ be a measured foliation on $S$. The {\it singular graph}  $K$ of $F$ is the union of the compact leaves that join the singular points of $F$. The {\it components} of $F$ are the (partial) measured foliations that are the closures of the connected components of $S\setminus K$. In this way, each measured foliation on $S$ can be decomposed as a union of finitely many components. The equivalence classes of the components of $F$ depend only on the equivalence class of $F$.

If $F$ and $G$ are two partial measured foliations on $S$ with disjoint supports, then their union can be naturally considered as a (partial) measured foliation on $S$, which we shall denote by $F+G$.
We shall say that $F$ is a {\it subfoliation} of $F+G$ and that the foliation $F+G$ (or any foliation equivalent to $F+G$) {\it contains} the foliation $F$ (or any foliation equivalent to $F$).

 \begin{lemma}\label{prop:int} Let $F$ and $G$ be two measured foliations on $S$. Then, the following are equivalent:
\begin{enumerate}
\item\label{prop:int1} $i(F,G)= 0$.
\item\label{prop:int2} $F\sim F'$ and $G\sim G'$, where $F'$ and $G'$ are partial measured foliations on $S$ such that $F'=F_1+F_2$ and $G'=G_1+G_2$ where $F_1$ and $G_1$ are equal as topological foliations, and where $F_2$ and $G_2$ have disjoint supports. (Some of the partial foliations $F_1$, $F_2$, $G_1$ and $G_2$ may be empty.)
\item \label{prop:int3}  $[F]$ is in the adherence set of $[G]$ in $\mathcal{UMF}$.
\end{enumerate}

\end{lemma} \begin{proof} The equivalence between (\ref{prop:int1}) and (\ref{prop:int2}) is well-known.
Let us prove that (\ref{prop:int2}) implies (\ref{prop:int3}).
Let $F'$ and $G'$ be partial measured foliations as in (\ref{prop:int2}). For the proof, we can assume that all four measured foliations, $F_1$, $F_2$, $G_1$ and $G_2$, are not empty. (In the contrary case, the proof is even simpler.) We then consider the measured foliation $F'+G_2$, and we let $[F'+G_2]$ be its equivalence class in $\mathcal{MF}$. For any sequence  $t_n$ of positive numbers convering to $0$, the sequences $[F'+t_nG_2]$ in $\mathcal{MF}$ converges to $[F']$. For the corresponding elements in $\mathcal{UMF}$, we have $[F'+t_nG_2]=[F'+G_2]$ for all $n$. This shows that in $\mathcal{UMF}$, $[F'+G_2]$ is in every neighborhood of $[F']$  (recall that a set in $\mathcal{UMF}$ is open if and only if its inverse image in $\mathcal{MF}$ is open). By the same argument, $[G'+F_2]$ is in every neighborhood of $[G']$ in $\mathcal{UMF}$. Since $[G'+F2]=[F'+G_2]$ in $\mathcal{UMF}$, $[F']$ is adherent to $[G']$. 
We now prove that (\ref{prop:int3}) implies (\ref{prop:int1}). The proof is by contradiction. Suppose that $i(F,G)\not= 0$. Then, by the continuity of the intersection function, we can find neighborhoods $N([F])$ of $[F]$ in $\mathcal{MF}$ and 
$N([G])$ of $[G]$ in $\mathcal{MF}$ such that $i(x,y)\not=0$ for all $x$ in $N([F])$ and for all $y$ in $N([G])$. We can furthermore suppose that $N([F])$ and $N([G])$ are saturated sts with respect to the equivalence relation on $\mathcal{MF}$ which identifies two equivalence classes of measured foliations whenever they can be represented by the same topological foliation. The images of $N([F])$ and of $N([G])$ in $\mathcal{UMF}$ are disjoint neighborhoods of the images of $F$ and of $G$ in that space. This shows that  $[F]$ is not in the adherence set of $[G]$ in $\mathcal{UMF}$.
\end{proof}

We already mentioned that Klarreich considered in her study of the Gromov boundary of the complex of curves $C(S)$ the subspace $\mathcal{UMF}'$ of $\mathcal{UMF}$ consisting of equivalence classes of minimal foliations. It is easy to see, by arguments analogous to those in the proof of Lemma \ref{prop:int}, that the space $\mathcal{UMF}'$  is Hausdorff.

We shall use the following

\begin{proposition}\label{prop:0} Let $F$ be a measured foliation and let $[F]$ beits image in $\mathcal{UMF}$. Then, the adherence set of $F$ is the set of equivalence classes in $\mathcal{UMF}$ of foliations $G$ which are of the form $G_1+G_2$ where $G_1$ is a sum of components of $F$ and where $G_2$ is a partial measured foliation whose support is disjoint from the support of a representative of $F$ by a partial measured foliation.
\end{proposition}

\begin{proof}
This is a direct consequence of the equivalence (\ref{prop:int1})$\Leftrightarrow$(\ref{prop:int2}) in Lemma \ref{prop:int}.
\end{proof}

A {\it system of curves} on $S$ is an isotopy class of a collection of pairwise non-homotopic disjoint essential simple closed curves on $S$.
Let $\mathcal{S}'$ be the set of isotopy classes of systems of curves on $S$.
There is a natural inclusion $j:\mathcal{S}' \hookrightarrow \mathcal{UMF}$ which is defined by associating to each element $C\in\mathcal{S}'$ the equivalence class of a partial measured foliation on $S$ whose support is the union of disjoint regular neighborhoods of the components of a system of curves on $S$ representing the  isotopy class $C$, and whose leaves are closed curves homotopic to these components. The choice of the total  transverse measure of each annulus does not matter.

We shall call a foliation on $S$ representing an element of  $\mathcal{UMF}$ which is the image of some element of $\mathcal{S}'$ by the map $j$ an {\it annular foliation}.

\begin{proposition}\label{prop:max} Let $F$ be a measured foliation on $S$ and let $[F]$ denote the corresponding element of $\mathcal{UMF}$. Then, 
$\mathcal{N}([F])=2^q -1$, where $q$ is the maximum, over all measured foliations $G$ containing $F$, of the number of components of $G$.
\end{proposition}

\begin{proof}
Let $G$ be a measured foliation containing $F$, let $k$ be the number of components of $G$, and let $[G_1],\ldots,[G_k]\in\mathcal{UMF}$ be the equivalence classes of these components.
Let $U(G)$ be the subset $\{[G_1],\ldots,[G_k]\}$ of $\mathcal{UMF}$ and let $V(G)$ be the set of all elements in $\mathcal{UMF}$ which represent foliations whose components are elements of $U(G)$. The cardinality of $V(G)$ is equal to the cardinality of the set of nonempty subsets of $U(G)$, which is $2^k-1$. 
 By Lemma \ref{prop:int}, any two elements in $V(G)$ are adherent. Thus, $\mathcal{N}([F])\geq 2^k -1$. It also follows from Lemma \ref{prop:int} that any adherence set in $\mathcal{UMF}$ that contains $F$ is a set of the form $V(G)$, for some measured foliation $G$ containing $F$. This completes the proof of the proposition.
 \end{proof}

Proposition \ref{prop:max} implies that the adherence number of any point in the space $\mathcal{UMF}$ is finite. Note that the maximum in the statement of that proposition is attained by any foliation $G$ obtained by taking a partial foliation $F'$ representing $F$ and completing it by an annular partial foliation supported in the complement of the support of $F'$ and containing the largest number possible of pairwise non-homotopic disjoint foliated annuli.

\begin{proposition}\label{prop:nb}
If $x\in\mathcal{UMF}$ is the class of an annular foliation, then $\mathcal{N}(x)= 2^q-1$, with $q=3g-3+p$. Furthermore, if
$y\in\mathcal{UMF}$ is not the class of an annular foliation, 
 then $\mathcal{N}(y)< \mathcal{N}(x)$.
 
\end{proposition}

\begin{proof}
This follows from Proposition \ref{prop:max}, and from the fact that if $F$ is annular, then the maximal number of components of a measured foliation containing it is $3g-3+p$, and if $F$ is not annular, the maximal number of components of a measured foliation $G$ containing it is $< 3g-3+p$.
\end{proof}

 The following is a consequence of Proposition \ref{prop:nb}, and it will be generalized in Proposition \ref{prop:Sk} below.
 
\begin{corollary}\label{co}
 Any homeomorphism of $\mathcal{UMF}$ preserves the image of $\mathcal{S}'$  in $\mathcal{UMF}$ by the map $j$. 
 \end{corollary}

\begin{proof}
A homeomorphism of $\mathcal{UMF}$ preserves adherence numbers of points.
\end{proof}

Let us now draw a another consequence of Proposition \ref{prop:nb} which will be useful in the proof of Proposition \ref{prop:Sk}  below.

\begin{corollary}\label{cor:not}
If $3g+p\not= 3g'+p'$, then $\mathcal{UMF}(S_{g,p})$ is not homeomorphic to $\mathcal{UMF}(S_{g',p'})$.
\end{corollary}

This corollary may also be obtained using a dimension count (one has to deal for that with dimension theory for non-Hausdorff spaces).

For every $k\geq 1$, we denote by of $\mathcal{S}_k$ the subset of $\mathcal{S}'$ consisting of homotopy classes that are representable by pairwise disjoint and non-homotopic $k$ simple closed curves. (In particular, $\mathcal{S}_1=\mathcal{S}$.)

\begin{proposition}\label{prop:Sk}
For any $k\geq 1$, any homeomorphism of $\mathcal{UMF}$ preserves the image of $\mathcal{S}_k$ by the map $j:\mathcal{S}'\to\mathcal{UMF}$.
\end{proposition}
\begin{proof}
Let $f$ be a homeomorphism of $\mathcal{UMF}$. By Corollary \ref{co}, $f$ preserves the subset $j(\mathcal{S}')$ of  $\mathcal{UMF}$. Note that $j(\mathcal{S}')$ is the disjoint union of the spaces $j(\mathcal{S}_k)$ with  $k=1,\ldots 3g-3+p$. It suffices to prove that if $m\not=n$ and for any $[F]\in j(\mathcal{S}_m)$ and $[G]\in j(\mathcal{S}_n)$, we have  $f([F])\not=[G]$. 
By Proposition \ref{prop:0}, the adherence set $\mathcal{A}([F])$ of $[F]$ (respectively $\mathcal{A}([G])$ of $[G]$) in $\mathcal{UMF}$ is homeomorphic to a finite union of spaces which are all homeomorphic to a space $\mathcal{UMF}(S_1)$ (respectively  $\mathcal{UMF}(S_2)$) where $S_1$ and  $S_2$ are (not necessarily connected) subsurfaces of $S$ 
which are the complement of the support of partial  foliations $F$  and $G$ respectively, representing $[F]$ and $[G]$ respectively. Since the number of components of $[F]$ and $[G[)$ are distinct, then if $(g_1,p_1)$ and $(g_2,p_2)$ are the genera and the number of punctures of $S_1$ and $S_2$, we have $3g_1+p_1\not=3g_2+p_2$.  By Corollary \ref{cor:not},  $\mathcal{A}([F])$ is not homeomorphic to $\mathcal([G])$. Thus, $f$ cannot send $[F]$ to $[G]$, which completes the proof of the proposition. 
\end{proof}

\section{Proof of Theorem \ref{th}}

We first  prove Statement (\ref{i1}) of Theorem \ref{th}.

Let $\mathcal{D}=j(\mathcal{S}')\subset\mathcal{UMF}$ and let $h$ be a homeomorphism of $\mathcal{UMF}$. Since $h$ preserves the set $j(\mathcal{S}_1)=j(\mathcal{S})$, $h$ induces a map from the vertex set $j(\mathcal{S})$ of $C(S)$ to itself. Furthermore, since, for each $k\geq 2$, $h$ preserves the set $j(\mathcal{S}_k)$ in $\mathcal{UMF}$, this action on the vertex set of $C(S)$ can be naturally extended to a simplicial automorphism $h'$ of $C(S)$.  Since  is not a sphere with at most four punctures, a torus with at most two punctures or a closed surface of genus two, it follows from the theorem of Ivanov, Korkmaz and Luo (see \cite{Ivanov}, \cite{Korkmaz}, \cite{Luo}) that the automorphism $h'$ is induced by an element $h''$ of the extended mapping class group $\Gamma^*(S)$. 
 The element $h''$ acts on $\mathcal{UMF}$, and, by construction, the action induced on $\mathcal{D}$ by this map is the same as that of $h$. This completes the proof of (\ref{i1}).
 
We now prove Statement (\ref{i2}) of Theorem \ref{th}, that is, if two elements of $\Gamma^*(S)$ have the same action on $\mathcal{D}$, then they are equal. It suffices to show that if
$g$ is an element of $\Gamma^*(S)$ that induces the identity map on  $\mathcal{D}$, then $g$ is the identity element. One can use the fact that the homomorphism from the extended mapping class group to the automorphism group of the complex of curves is injective, but we can give a direct argument as follows.

Suppose that $g$ induces the identity map on $\mathcal{D}$. Since the restriction of the quotient map $\mathcal{PMF}\to\mathcal{UMF}$ to the natural  image of $\mathcal{S}'$ in $\mathcal{PMF}$ is injective, the action of $g$ on $\mathcal{PMF}$ induces the identity map on the natural image of $\mathcal{S}'$ in $\mathcal{PMF}$. We show that this implies that $g$ is the identity element of $\Gamma^*(S)$. This will follow directly from Thurston's classification of mapping classes.

First, suppose that $g$ is orientation-preserving. Then, Thurston's classification of the elements of the mapping class group also gives a description of the dynamics of $g$ on the space $\mathcal{PMF}$. If $g$ is not of finite order, then, since $S$ is not a sphere with at most foir punctures or a torus with at most two punctures, the fixed point set of $g$ on $\mathcal{PMF}$ has codimension $>1$.  Thus, $g$ cannot fix pointwise the image of $\mathcal{S}$ in $\mathcal{PMF}$. Therefore, $g$ is of finite order. Suppose that $g$ is not the identity element. Then, we can find a homeomorphism $\gamma$ of $S$ which has finite order and which represents $g$. The fixed point set of $g$ in $\mathcal{PMF}$ can be studied by examining the projective measured foliation space of the quotient surface $S/\gamma$ which embeds as a set of codimension $>1$ in $\mathcal{PMF}(S)$. Since $h$ is not a hyperelliptic involution of the closed surface of genus 2, the fixed point set of $g$ in $\mathcal{PMF}(S)$ has also codimension $>1$. Therefore $g$ cannot fix every element in the image of $\mathcal{S}'$ in $\mathcal{PMF}$. Thus, $g$ is the identity element of $\Gamma^*(S)$.

Finally, suppose that $g$ is orientation-reversing.
If $g$ induces the identity on $\mathcal{D}$, the same is true for $g^2$, which is orientation-preserving, and by the previous discussion, $g^2$ is the identity. Therefore $g$ is an involution. The fixed point set of the action of an involution in $\Gamma^*(S)$ on $\mathcal{PMF}$ can be studied by looking at the projective measured foliation of the quotient surface of $S$ by a homeomorphism of order two, and this fixed point set has codimension $>1$ provided $g$ is not the identity. Therefore, if $g$ fixes every point in $\mathcal{D}$, it is the identity. This completes the proof of the theorem.

\end{document}